\documentclass[11pt]{article}

\usepackage{graphicx}
\usepackage{textcomp}
\usepackage{amsmath}
\usepackage{amsfonts}
\usepackage{epsfig}

\newcommand{\be}{\begin{equation}}
\newcommand{\ee}{\end{equation}}
\newcommand{\beqn}{\begin{eqnarray}}
\newcommand{\eeqn}{\end{eqnarray}}
\newcommand{\beqns}{\begin{eqnarray*}}
\newcommand{\eeqns}{\end{eqnarray*}}
\newcommand{\lkr}{\left(}
\newcommand{\lkv}{\left[}
\newcommand{\rkv}{\right]}
\newcommand{\rkr}{\right)}
\newcommand{\lfi}{\left\{}
\newcommand{\rfi}{\right\}}

\newcommand{\fr}[1]{(\ref{#1})}

\newcommand{\ro}{\varrho}
\newcommand{\ph}{\varphi}

\newcommand{\af}{\alpha}
\newcommand{\eps}{\varepsilon}

\newcommand{\ga}{\gamma}

\newcommand{\om}{\omega}
\newcommand{\lam}{\lambda}

\newcommand{\Om}{\Omega}

\newcommand{\ubeta}{\underline{\beta}}
\newcommand{\uaf}{\underline{\alpha}}
\newcommand{\unp}{\underline{p}}
\newcommand{\unu}{\underline{u}}
\newcommand{\barq}{\bar{q}}

\newcommand{\hf}{\hat{f}_n}

\newcommand{\psimjk}{\psi_{mjk}}
\newcommand{\phmlk}{\ph_{mj_0k}}
\newcommand{\psijk}{\psi_{jk}}

\newcommand{\fm}{f_m}

\newcommand{\ym}{y_{m}}
\newcommand{\gm}{g_{m}}

\newcommand{\gmu}{g_{m}(u)}

\newcommand{\ajk}{a_{j_0k}}
\newcommand{\bjk}{b_{jk}}
\newcommand{\hajk}{\widehat{a}_{j_0k}}
\newcommand{\hbjk}{\widehat{b}_{jk}}

\newcommand{\hfm}{\widehat{f}_m }

\newtheorem{theorem}{Theorem}
\newtheorem{lemma}{Lemma}

\newcommand{\Bpqsa}{B_{p,q}^s (A) }

\newcommand{\sumk}{\sum_{k=0}^{2^j-1}}

\newcommand{\sumM}{\sum_{l=1}^M}

\newcommand{\jo}{{j_0}}

\newcommand{\sstar}{{s^*}}

\newcommand{\EE}{\ensuremath{{\mathbb E}}}
\newcommand{\II}{\ensuremath{{\mathbb I}}}
\newcommand{\PP}{\ensuremath{{\mathbb P}}}

\newcommand{\ints}{\ensuremath{{\mathbb Z}}}

\newcommand{\RR}{\ensuremath{{\mathbb R}}}

\begin{document}

\title{\bf Multichannel Boxcar Deconvolution with Growing Number of Channels}

\author{{\em Marianna ~Pensky},\\
         Department of Mathematics,
         University of Central Florida,\\
         Orlando, FL 32816-1364, USA.\\
          \\
        {\em Theofanis ~Sapatinas},\\
        Department of Mathematics and Statistics,
        University of Cyprus,\\
        P.O. Box 20537,
        CY 1678 Nicosia,
        Cyprus.}

\date{}

\bibliographystyle{plain}
\maketitle

\begin{abstract}

We consider the problem of estimating the unknown response function in the  multichannel deconvolution model with  a boxcar-like kernel which is of particular interest in signal processing.  It is  known that, when the number  of channels is finite, the precision of reconstruction of the response function increases as the number of channels $M$ grow (even when the total number of observations $n$ for all channels $M$ remains constant) and this requires that the parameter of the channels form a Badly Approximable  $M$-tuple.

Recent advances in data collection and recording techniques made it of urgent interest  to study the case 
when the number of channels $M=M_n$ grow with the total number of observations $n$. However, in real-life situations, the number of channels $M = M_n$ usually refers to the number of physical devices 
and, consequently,  may grow to infinity only at a slow  rate as $n \rightarrow \infty$. Unfortunately, existing theoretical results  cannot be blindly applied to accommodate the case when $M=M_n \rightarrow \infty$ as $n \rightarrow \infty$.  This is due to the fact that, to the best of our knowledge, so far no one have studied the construction of a Badly Approximable $M$-tuple of a growing length on a specified interval, of a non-asymptotic length, of the real line, as $M$ is growing. Therefore,  this generalization requires  non-trivial results in number theory.

When  $M=M_n$ grows slowly as  $n$ increases, we develop a procedure for the construction of a Badly Approximable $M$-tuple on a specified interval, of a non-asymptotic length, together with a lower bound associated with this $M$-tuple, which explicitly shows its dependence on $M$ as $M$ is growing. This result is further used for the evaluation of the $L^2$-risk of the suggested adaptive wavelet thresholding estimator of the unknown response function and, furthermore, for the choice of the optimal number of channels $M$ which minimizes the $L^2$-risk. 


\vspace{3mm} {\bf  Keywords}: {Adaptivity, badly approximable tuples, Besov spaces, 
Diophantine approximation, functional deconvolution, Fourier analysis, Meyer wavelets, nonparametric estimation,  wavelet analysis}

\vspace{3mm}{\bf AMS (2000) Subject Classification}: {Primary: 62G05, 11K60, Secondary: 62G08, 35L05}
\end{abstract}

\section {Introduction}
\label{sec:intro}
\setcounter{equation}{0}

We consider the estimation problem of the unknown response function
$f(\cdot) \in L^2(T)$  from observations $y(u_l, t_i)$, $l=1,2,\ldots,M$, $i=1,2,\ldots,N$, where
\be
y(u_l, t_i) = \int_T g(u_l, t_i-x)f(x) \, dx + \eps_{li},\ \
\  u_l \in U, \;\; t_i = (i-1)/N,
\label{convdis}
\ee
where $U=[a,b]$, $0 < a < b < \infty$, $T=[0,1]$ and  $\eps_{li}$ are standard Gaussian random variables, 
independent for different $l$ and $i$.
We shall  be interested in the case when the blurring (or kernel) function $g(\cdot,\cdot)$ is the, so called, boxcar-like kernel, i.e.,
$$
g(u,t) =\frac{ \gamma(u)}{2}\, \II(|t| < u),
$$
where $\gamma(\cdot)$ is some positive function such that  
\begin{equation} \label{bc-gammacond} 
\gamma_1 \leq \gamma(u) \leq \gamma_2,\quad u\in U,
\end{equation}
for some $0 < \gamma_1 \leq \gamma_2 < \infty$. (Obviously, this is
true if $\gamma(\cdot)$ is a continuous function.) 
Hence, \fr{convdis} is of the form 
\begin{equation}  \label{eq:box-car1}
y(u_l, t_i) = \frac{\gamma(u_l)}{2}\ \int_0^1  \II(|t_i-x| < u_l) f(x)\,  dx + \eps_{li},\ \
\  u_l \in U, \;\; t_i = (i-1)/N,
\end{equation}
for $l=1,2,\ldots,M$ and $i=1,2,\ldots,N$.

In signal processing, this model is referred to as a multichannel deconvolution model,
where $M$ is the number of channels and $N$ is the number of observations per channel, so that
$n = MN$ is the total number of observations. 
We assume that the measurements  $t_i$ in each channel are equispaced but  
the observer can choose the number of channels $M$ and the points $u_l$, $l=1, \ldots, M$, in \fr{convdis}  
prior to the experiment as a part of experimental design.
In order to be able to access convergence rates depending on the number of channels $M$ and the choice of 
points $u_l$, $l=1,2, \ldots, M$, we shall further assume that the total number of observations $n$ is fixed and very large (i.e., $n \rightarrow \infty$). 
The objective is to choose  $M$ and $u_l$, $l=1,2, \ldots, M$, which ensure the construction of an estimator of the response function $f$ 
with the highest possible convergence rates in terms of $n$.

Note that  standard deconvolution (i.e., when $a=b$) with the boxcar kernel  
(i.e., when $\ga(u) =1/u$, for some {\em fixed} $u >0$)
is a common model in many areas of
signal and image processing which include, for instance, LIDAR
remote sensing and reconstruction of
blurred images. LIDAR is a lazer device which emits pulses,
reflections of which are gathered by a telescope aligned with the
lazer, see, e.g., Park {\em et al.} (1997) and Harsdorf
\& Reuter (2000). The return signal is used to determine distance
and the position of the reflecting material. However, if the system
response function of the LIDAR is longer than the time resolution
interval, then the measured LIDAR signal is blurred and the
effective accuracy of the LIDAR decreases. This loss of precision can be corrected by deconvolution. In
practice, measured LIDAR signals are corrupted by additional noise
which renders direct deconvolution impossible. If $M\geq2$,  then we talk
about a multichannel deconvolution model with blurring functions $g_l(t) = g(u_l, t)$.

Although standard deconvolution models are traditionally solved using the Fourier transform
or the Fourier series,  if the corresponding blurring function $g(\cdot)$ is a boxcar-like kernel, implementation of the standard Fourier series based technique is impossible.
This happens when the Fourier transform of $g(\cdot)$ has real
zeros, e.g., when $g(\cdot)$ is the boxcar kernel $g(x) = (2u)^{-1} \II(|x| \leq u)$, for some {\em fixed} $u>0$.  When $M=1$,  Johnstone {\em et al.} (2004) and Johnstone \& Raimondo (2004) 
managed to circumvent this obstacle by considering a boxcar kernel $g(\cdot)$ with irrational scale.   
Their method is based on the fact that the Fourier coefficients of the boxcar kernel do not vanish
at frequencies $(\pi k u)$ when $u$ is a {\em Badly Approximable} (BA) number.
An irrational number $u$ is BA if the terms $a_n=a_n(u)$ of its continued fraction expansion $[a_0;a_1,a_2 \ldots]$, where $a_0$ is an integer and $a_1,a_2,\ldots$ is an infinite sequence of positive integers, are bounded, i.e., $\sup_n a_n(u) < \infty$. This notion is related to the fact that a BA number cannot be approximated well by a rational number which leads 
to the fact that $f$ can be recovered reasonably well. Since standard deconvolution is a particular example of linear statistical ill-posed inverse problems in the sense of Hadamard, i.e., the inversion does not depend continuously on the observed data, Johnstone \& Raimondo (2004) used number theory to prove that  the degree of ill-posedness in boxcar deconvolution  is $\nu=3/2$.
Roughly speaking, the degree of ill-posedness specifies how much the error in the right-hand side of the equation is amplified in the solution.
For example, if $f$ belongs to a space with a smoothness index $s>0$ and the degree of ill-posedness is $\nu>0$, then, the quadratic risk of the best possible estimator of the response function $f$   is of the order $O \lkr n^{-\frac{2s}{2s + 2\nu +1}} \rkr$.

De Canditiis \& Pensky (2004, 2006), following mathematical ideas of Casey \& Walnut (1994),
extended the results of Johnstone {\em et al.} (2004) and Johnstone \& Raimondo (2004), and showed that if $M$ is finite, $M \geq 2$, one of the $u_l$'s is a BA number, and $u_1,u_2, \ldots, u_M$ is a BA
 $M$-tuple, then the degree of ill-posedness is  $\nu = 1 + 1/(2M)$. 
The notion of a BA  $M$-tuple refers to a collection of $M$ irrational numbers 
which are difficult to approximate simultaneosly by fractions with the same denominator. 
It will be discussed in depth in Section \ref{sec:background}.  
Therefore, in the case of $M$ channels, the estimation problem requires a construction of a BA $M$-tuple which has been 
accomplished by the number theory community (it is described in, e.g., Schmidt (1969, 1980)).

Recent advances in data collection and recording techniques made it of urgent interest  to study the case 
when the number of channels $M=M_n$ grow with the total number of observations $n$. 
It turns out that when the number of channels $M=M_n$ grows fast as the total number of observations $n$ increases, 
one does not need to make a special choice of the points  $u_l$, $l=1,2 \ldots, M$, and it is sufficient to take them to be equidistant. Indeed,  Pensky \& Sapatinas (2010) considered the discrete multichannel deconvolution model \fr{convdis} as observations  on the continuous functional 
deconvolution model 
\be
y(u, t) = f*g (u, t) + \frac{1}{\sqrt{n}}\,
z(u,t),\ \ \  u \in U, \;\; t \in T,
\label{convcont}
\ee
where
$z(u,t)$ is assumed to be a two-dimensional Gaussian white noise, i.e., a
generalized two-dimensional Gaussian field with covariance function
 $ 
\EE [z (u_1, t_1) z (u_2, t_2)] = \delta(u_1-u_2) \delta(t_1-t_2),
$ 
where $\delta(\cdot)$ denotes the Dirac $\delta$-function, and 
$ 
f*g(u, t) = \int_T g(u,t-x)f(x)\, dx 
$ 
with the blurring (or kernel) function $g(\cdot,\cdot)$ assumed to be known.
 If $a=b$, the functional deconvolution model  (\ref{convcont}) reduces to the standard
deconvolution model which attracted attention of a number of
researchers, e.g., Donoho (1995), Abramovich \& Silverman (1998), 
Kalifa \& Mallat (2003), Johnstone\ {\em et al.} (2004), Donoho \& Raimondo (2004), Johnstone  \& Raimondo
(2004), Neelamani\ {\em et al.} (2004),   Kerkyacharian\ {\em et al.} (2007), 
Cavalier  \& Raimondo  (2007) and Chesneau (2008), among others.

Formulation of the functional deconvolution model (\ref{convcont}) allowed  Pensky \& Sapatinas (2010)
to study the interplay between discrete and continuous deconvolution models. 
The {\em ideal} continuous deconvolution model (\ref{convcont}) assumes that one can measure $y(u,t)$ at 
{\em any}\, $u \in U$ and $t \in T$ and that $n^{-1/2}$ marks the precision of these observations. 
Nevertheless, this  does not happen in real-life situations where one observes $y(u,t)$ {\em only}\, at  
the points $ u_l \in U, \;\; t_i = (i-1)/N $ for $l=1,2,\ldots,M$ and $i=1,2,\ldots,N$.
Pensky \& Sapatinas (2010) showed that the degree of ill-posedness in the continuous deconvolution model \fr{convcont} is $\nu=1$ 
and that it can be attained in the discrete deconvolution model \fr{convdis} if $M = M_n \geq c_0 n^{1/3}$ for some constant $c_0 >0$, independent of $n$.
Indeed, in this case, one does not need to employ BA numbers or BA $M$-tuples: it is sufficient to observe the discrete deconvolution model  \fr{convdis} 
at equidistant points $u_l = a+ l (b-a)/M$, $l=1,2, \ldots, M$. This set up provides the ``best possible'' minimax convergence rates (under the $L^2$-risk and over a wide range of Besov balls) in the model.

However, in real-life situations, the number of channels $M$ usually refers to the number of physical devices and, consequently, cannot be very big. Therefore, $M = M_n \geq c_0 n^{1/3}$ may be impossible although it is natural to assume that $M=M_n$ may grow to infinity at a slower rate as $n \rightarrow \infty$.  Unfortunately, the theoretical results obtained  by De Canditiis \& Pensky (2006)
cannot be blindly applied to accommodate the case when $M=M_n \rightarrow \infty$ as $n \rightarrow \infty$.  This is due to the fact that, to the best of our knowledge, so far no one have studied the construction of a BA  $M$-tuple of a growing length on a specified interval, of a non-asymptotic length, of the real line. Therefore, this generalization requires  non-trivial results in number theory.

Our aim is to investigate the situation when $M=M_n$ grows slowly with $n$ 
and to derive necessary new results  in number theory in order to devise a technique which allows to 
approach minimax convergence rates (under the $L^2$-risk and over a wide range of Besov balls) in the continuous model \fr{convcont} with a factor which grows slower than any power of $n$.
This situation seems to be of a particular interest nowadays  since data recording equipment is 
getting cheaper and cheaper while overall volumes of data is growing very fast.

When  $M=M_n$ grows slowly as  $n$ increases, we develop a procedure for the construction of a BA $M$-tuple on a specified interval, of a non-asymptotic length,   together with a lower bound associated with this $M$-tuple, 
which explicitly shows its dependence on $M$ as $M$ is growing. This result is further used for evaluation of the $L^2$-risk 
of the suggested adaptive wavelet thresholding estimator of the unknown response function and, furthermore, for the choice of the optimal number of channels $M$ which minimizes the $L^2$-risk. 
 
 The theoretical results that we have obtained provide a cross-area between number theory, statistics and signal processing. We hope to alert the number theory community to a new problem of constructing a BA   $M$-tuple on a specified interval, of a non-asymptotic length, of the real line, as $M$ is growing. 
On the other hand, we believe that our findings will also be of interest to researchers in statistics and signal processing.

The rest of the paper is organized as follows. Section \ref{sec:background} provides some number theory background which is required for understanding the material presented in subsequent sections. Section \ref{sec:est_method} briefly reviews the adaptive wavelet thresholding estimator 
introduced in Pensky \& Sapatinas (2009). Section \ref{sec:relation} explains the relationship between the $L^2$-risk of the estimator 
obtained in Section \ref{sec:est_method} and the theory of Diophantine approximation, thus, motivating the derivation of the
new results in number theory obtained in Section \ref{sec:num_theory}. In particular, the objective of Section \ref{sec:num_theory} 
is the construction of a BA  $M$-tuple on a specified interval when $M=M_n \rightarrow \infty$  as $n \rightarrow 0$ and the development of 
related asymptotic bounds which are necessary in order to choose an optimal value of $M=M_n$ in this case. 
Section \ref{sec:asymp_bounds} provides the asymptotic upper bounds for the $L^2$-risk of the adaptive wavelet thresholding estimator constructed in Section
\ref{sec:est_method} when $M=M_n$ is a slowly growing function of $n$. We conclude in Section  \ref{sec:discussion} with a brief discussion while Section \ref{sec:proofs} contains the proofs of the theoretical results obtained in earlier sections.

\section {Background Results in Number Theory}
\label{sec:background}
\setcounter{equation}{0}

The theory of Diophantine approximation is an important branch of number theory
(see, e.g., Edixhoven \& Evertse  (1993),  Lang  (1966), Masser {\it et al.}  (2003) and Schmidt (1980, 1991)). 
One important topic of the above theory is the simultaneous approximation of linear forms, which was pursued as early as 
mid-19th century by Dirichlet and later studied by a number of profound researchers in the field.
In particular, it is known that for any real numbers $\beta_1 , \beta_2, \ldots,  \beta_M$ there exist 
integer numbers $q$ and $p_1, p_2, \ldots, p_M$ such that 
\begin{equation}
\label{eq:fan-orig}
\max_{i=1,2,\ldots, M}  | \beta_i q - p_i | < \frac{M}{(M+1)} |q|^{-1/M}.
\end{equation}
The above result was proved by Minkowski  and has been expanded in the recent 
years to cover systems of linear forms (see, e.g., Schmidt (1980), p. 36, pp. 40-41). We note that in the case where $M=1$, 
the constant $C(M)=M/(M+1)$ in (\ref{eq:fan-orig}) reduces to 1/2 whereas, by Hurwitz's theorem, the best possible value is 
$1/\sqrt{5}$ (see, e.g., Schmidt (1980), Theorem 2F, p. 6). For $M=2$, $C(M)$ takes the value 2/3; the best possible value is unknown
although if $C_0(2)$ denotes the infimum of admissible values of $C(M)$ for $M=2$, then it is known that $\sqrt{2/7} \leq C_0(2) \leq 0.615$
(see, e.g., Schmidt (1980), p. 41). Furthermore, the corresponding best constant in the case of systems of linear forms is positive, meaning that 
it cannot be replaced by arbitrary small constants (see, e.g., Schmidt (1980), Section 4, pp. 41-47) .

We, however, are interested in the opposite result. Namely, the real numbers $\beta_1, \beta_2, \ldots,  \beta_M$
form a {\it BA $M$-tuple} if for any integer numbers  $q>0$ and $p_1, p_2, \ldots, p_M$
one has 
\be  \label{BAtuple}
\max_{i=1,2,\ldots, M} | \beta_i q - p_i | \geq B(M) q^{-1/M}, 
\ee
for some constant $B(M)>0$,  dependent on $M$ (and $\beta_1, \beta_2,\ldots,\beta_M$) but independent of $q$ and $p_1, p_2, \ldots, p_M$ 
(see, e.g., Schmidt (1980), p. 42).  It is well-known  that the set of all BA $M$-tuples has Lebesgue measure zero, 
but nevertheless this set is quite large, namely there are uncountably many BA $M$-tuples (see Cassels (1955), Davenport (1962)) and the 
Hausdorff dimension of the set of all BA $M$-tuples is equal to $M$ (see Schmidt (1969)). 
In the case where $M=1$, the number  $\beta=\beta_1$ which satisfies (\ref{BAtuple}) is referred to in the Diophantine 
approximation literature as a BA number (see, e.g., Schmidt (1980), p. 22); in view of Hurwitz's theorem, the constant $B(1)$ in this case must satisfy $0 < B(1)<1/\sqrt{5}$ (see, e.g., Schmidt (1980), pp. 41-42). Furthermore, a characterization result  exists, namely a real number, that is not an integer, is BA if and only if its continued fraction coefficients are bounded. The latter is often used as a definition of a BA number, 
however, there is no analogous characterization  for $M>1$ (see, e.g., Schmidt (1980), Theorem 5F, p. 22). The above definitions of BA numbers and BA $M$-tuples have been also extended to cover BA systems of linear forms (see, e.g., Schmidt (1980), pp. 41) and their existence was proved by Perron, providing also an algorithm for constructing BA linear forms (see, e.g., Schmidt (1980), Theorem 4B, p. 43). Furthermore, it has been established the existence of uncountably many BA systems of linear forms (see Schmidt (1969)).

In what follows, we are interested in the case of BA $M$-tuples. Although, as indicated above, an algorithm is available for constructing BA $M$-tuples on the real line, these do not necessarily belong to any specified interval of the real line. Furthermore,  if $M$ is strictly fixed (independent of $q$), one can treat $B$ in (\ref{BAtuple}) as a positive constant; this, however, becomes impossible if the value of $M$ is growing. Using the technique described in Schmidt (1980), Section 4, pp. 43-45, we show that one can construct a BA $M$-tuple
$\beta_1, \beta_2, \ldots,  \beta_M$ of real numbers so that it lies in any specified interval $(a,b)$, $a < b$, of non-asymptotic length, of the real line, and derive a lower bound for $B(M)$ in (\ref{BAtuple}) as $M \rightarrow \infty$. 
This result is proved in Section \ref{sec:num_theory}.

\section {An adaptive  wavelet thresholding estimator}
\label{sec:est_method}
\setcounter{equation}{0}

Let  $\ph^*(\cdot)$ and $\psi^*(\cdot)$ be the Meyer scaling and
mother wavelet functions, respectively, in the real line  (see, e.g., Meyer (1992) or
Mallat (1999)). As usual,  
$$
\ph^*_{jk}(x) = 2^{j/2}\ph^*(2^jx-k), \quad  \psi^*_{jk}(x) =
2^{j/2}\psi^*(2^jx-k), \quad j,k \in \ints, \quad x \in \RR,
$$
are, respectively, the dilated and translated Meyer scaling and
wavelet (orthonormal) basis functions at resolution level $j$ and scale position
$k/2^j$.  Similarly to Section 2.3 in Johnstone  {\it et al.} (2004), 
we obtain a periodized version of the Meyer
wavelet basis, by periodizing the basis functions
$\{\ph^*(\cdot),\psi^*(\cdot)\}$, i.e., for $j \geq 0$ and $k=0,1,\ldots,2^j-1$,
$$
\ph_{jk}(x) = \sum_{i \in \ints} 2^{j/2} \ph^*(2^j (x +i) -  k),
\quad \psi_{jk}(x) = \sum_{i \in \ints} 2^{j/2} \psi^*(2^j (x +i) -
k), \quad x \in T.
$$

Let $\langle \cdot,\cdot \rangle$ denote the inner
product in the Hilbert space $L^2(T)$ (the space of
squared-integrable functions defined on $T$),
i.e., $\langle f,g \rangle = \int_{T}f(t)\overline{g(t)} dt$ for
$f,g \in L^2(T)$. Let $e_m(t) = e^{i 2 \pi m t}$, $m \in \ints$, and let $f_m =  \langle e_m, f  \rangle$, $\gmu = \langle e_m, g(u,\cdot) \rangle$, $u \in U$.
For any $j_0 \geq 0$ and any $j \geq j_0$, let
$\phmlk  =   \langle e_m, \ph_{j_0k} \rangle$ and $ \psimjk = \langle
e_m, \psijk \rangle$,  where $\{\phi_{j_0,k}(\cdot),\psi_{j,k}(\cdot)\}$  is
the periodic Meyer wavelet basis  introduced above. 

Using the periodized Meyer wavelet basis described above,
and for any $j_0 \geq 0$, any (periodic) $f(\cdot) \in L^2(T)$ can
be expanded as 
\be 
f(t) = \sum_{k=0}^{2^{j_0}-1} a_{j_0k} \ph_{j_0k}
(t) + \sum_{j=j_0}^\infty \sum_{k=0}^{2^j -1} b_{jk} \psijk (t), \quad t \in T.
\label{funf} 
\ee 
Furthermore, by Plancherel's formula, the scaling
coefficients, $\ajk=\langle f, \ph_{j_0k} \rangle$, and the wavelet
coefficients, $\bjk=\langle f,\psi_{jk}\rangle$, of $f(\cdot)$ can
be represented as 
\be 
\ajk = \sum_{m \in C_{j_0}} \fm \overline{\phmlk},\ \ \ \bjk = \sum_{m \in C_j} \fm
\overline{\psimjk}, 
\label{alkandblk} 
\ee 
where $C_{j_0} = \lfi m:
\phmlk \neq 0 \rfi$ and, for any $j \geq j_0$, $C_j = \lfi m:
\psimjk \neq 0 \rfi$. Note that both $C_{j_0}$ and $C_{j}$, $j \geq j_0$, are
subsets of $(2\pi/3) [-2^{j+2}, -2^j] \cup [2^j, 2^{j+2}]$, i.e.,
\be  \label{range_of_m}
 |m| \in  (2\pi/3)\ [2^j, 2^{j+2}]
\ee
due to the fact that Meyer wavelets are band
limited (see, e.g., Johnstone {\it et al.} (2004), Section 3.1).

Reconstruct   the unknown   response function $f(\cdot) \in L^2(T)$ in \fr{eq:box-car1}  as
\be \label{fest}
\hat{f}_n (t) = \sum_{k=0}^{2^{j_0} -1} \hajk
\ph_{j_0k} (t) + \sum_{j=j_0}^{J-1} \sum_{k=0}^{2^j-1}\hbjk \II(|\hbjk|
\geq \lam_{j}) \psijk (t), \quad t \in T,  
\ee 
where $\hajk$ and $\hbjk$ are the natural estimates of $\ajk$ and $\bjk$, respectively (see \fr{funf} and \fr{alkandblk}), given by
\begin{equation} 
\hajk = \sum_{m \in C_{j_0}} \hfm
\overline{\phmlk},\ \ \ \hbjk = \sum_{m \in C_j} \hfm
\overline{\psimjk}. 
\label{coefest} 
\end{equation}
with $\hfm$ obtained by
\begin{equation*} \hfm =
\bigg ( \sumM \overline{\gm (u_l)} \ym(u_l) \bigg ) \Big/ \bigg (
\sumM |\gm(u_l)|^2 \bigg ), \quad u_l \in U, \ \  l=1,2,\ldots,M;
\label{fmexprd} 
\end{equation*}
here $g_m(u_l)$ and $y_m(u_l)$, $l=1,2,\ldots,M$, are the discrete Fourier coefficients of $y(u,\cdot)$ and $g(u,\cdot)$, respectively,
obtained by applying the discrete Fourier transform to the equation \fr{eq:box-car1}.
Note that, in this case, 
\begin{equation}  \label{gmul}
g_0(u_l)=1 \quad \text{and} \quad g_m(u_l) = \gamma(u_l) \frac{\sin(2 \pi m u_l)}{2\pi  m }, \ \ m \in \ints \setminus \{0\}\ \ l=1,2,\ldots, M.
\end{equation}

The choices of the resolution levels $j_0$ and $J$ and the thresholds $\lam_{j}$ will be described in Section \ref{sec:asymp_bounds}  when we examine an
expression for the $L^2$-risk of the estimator \fr{fest} over a collection of Besov balls, leading to an adaptive estimator (i.e., its construction is independent of the Besov ball parameters that are usually unknown in practice).

Among the various characterizations of Besov spaces for periodic
functions defined on $L^p(T)$ in terms of wavelet bases, we recall
that for an $r$-regular ($0 < r \leq \infty$) multiresolution analysis with $0< s < r$ and
for a Besov ball $B_{p,q}^s (A)$ of radius $A>0$ with $1 \leq p,q
\leq \infty$, one has that, with $s' = s+1/2-1/p$,
\begin{equation}
\label{bpqs}
B_{p,q}^s (A) = \Bigg \{ f(\cdot) \in L^p(T): \Bigg (
\sum_{k=0}^{2^{j_0}-1}|a_{j_{0}k}|^p \Bigg )^{\frac{1}{p}} 
+ \Bigg ( \sum_{j=j_0}^{\infty} 2^{js'q} \bigg ( \sumk
|\bjk|^p \bigg )^\frac{q}{p} \Bigg )^{\frac{1}{q}} \leq A \Bigg \},
\end{equation}
with respective sum(s) replaced by maximum if $p=\infty$ or
$q=\infty$ (see, e.g., Johnstone {\it et al.} (2004), Section 2.4). (Note that, for the Meyer wavelet basis,
considered considered above, $r=\infty$.)

The parameter $s$ measures the number of derivatives, where the existence
of derivatives is required in an $L^{p}$-sense, while the parameter $q$
provides a further finer gradation. The Besov spaces include, in particular, the well-known Sobolev and 
H\"{o}lder spaces of smooth functions 
but in addition less traditional spaces, like the space of functions
of bounded variation.
The latter functions are of statistical interest because they allow for 
better models of spatial inhomogeneity 
(see, e.g., Meyer (1992)).

The precision of the estimator \fr{fest} is measured by the (maximal) $L^2$-risk given by 
\be
\label{target}
R_n (\hf) = \sup_{f \in B_{p,q}^s (A)} \EE \|\hf -f\|^2_2. 
\ee
We are interested in the asymptotic rate of convergence of the estimator $\hf$, i.e., we are interested in the following asymptotical upper bounds  
$$R_n (\hf)\leq C \gamma_n \quad \text{as} \quad  n \rightarrow \infty,$$
where $\{\gamma_n\}_{n=1}^{\infty}$ is a positive sequence converging to 0 as $n \rightarrow \infty$
and $C>0$ is a generic constant, independent of $n$, which may take different values at
different places.

Hereafter,  $\| \cdot \|_2$ denotes the $L^2$-norm,  $\hat{f}_n(\cdot)$ is an estimator 
(i.e., a measurable function) of $f(\cdot) \in L^2(T)$, based on observations from model \fr{eq:box-car1}, 
and the expectation in \fr{target} is taken under the true $f(\cdot)$.

\section{Relation to the theory of Diophantine approximation}
\label{sec:relation}
\setcounter{equation}{0}

By direct evaluations (see also Lemma \ref{lemm:bound1} and its proof in the Appendix),
one can show that
$$
\EE |\hbjk - \bjk|^2 =  N^{-1} \sum_{m \in C_j} |\psimjk|^2\ \Big [ \sumM
|\gm(u_l)|^2 \Big ]^{-1}.
$$
Since in the case of Meyer wavelets, $|\psimjk| \leq 2^{-j/2}$ and $|C_j| \asymp 2^j$ 
(see, e.g., Johnstone {\it et al} (2004), p. 565), we derive  
\begin{equation} \label{power2}
\EE |\hbjk - \bjk|^2 =   O \Big ( n^{-1} \Delta_1 (j) \Big ),
\end{equation}
where
$$
\Delta_1 (j) = \frac{1}{|C_j|} \sum_{m \in C_j} \lkv M^{-1} \sum_{l=1}^M |g_m(u_l)|^2 \rkv^{-1} \equiv
\frac{1}{|C_j|} \sum_{m \in C_j} [\tau_1 (m)]^{-1}
$$
with $\tau_1 (m) \equiv \tau_1 (m; \underline{u},M)   = M^{-1} \sum_{l=1}^M |g_m(u_l)|^2$, $\unu = (u_1, u_2,\ldots, u_M)$.
By \fr{bc-gammacond} and \fr{gmul}, one has
\be \label{eq:box-car-taulower}
\tau_1 (m)   \asymp \frac{1}{m^2\, M}   \sum_{l=1}^M \sin^2(2 \pi m u_l).
\ee
Here,  $u(m) \asymp v(m)$ mean that there exist constants $C_1>0$ and $C_2>0$,
independent of $m$, such that $0<C_1 v(m) \leq u(m) \leq C_2 v(m)<\infty$ for every $m$.

Therefore, the risk of the estimator $\hat{f}_n (t)$ defined in  \fr{fest}
is determined  by the rate of growth of $\Delta_1 (j)$ as $j \rightarrow \infty$ 
which, in turn, depends on the rate at which $\tau_1 (m)$ goes to zero as $m \rightarrow \infty$.

It is easy to see that for some choices of $M$ and $\unu$ (e.g., $M=1$, $u_1=u=1$), one has 
$\min_m \tau_1 (m; \underline{u},M) =0$ for every $m$ which leads to an infinite variances of the estimated coefficients 
$\hbjk$ and, consequently, to an infinite   $L^2$-risk.
Hence,  the choice of $M$ and the selection of points $\unu$ is of an uttermost importance.
In particular, we want to choose points $(u_1, u_2,\ldots, u_M)$ such that
$\sum_{l=1}^M \sin^2(2 \pi m u_l)$ is as large as possible for $m \in C_j$ and large $j$.

Moreover, for any choice of $M$ and any selection of points $\underline{u}$, one has
$\tau_1 (m; \underline{u},M) \leq K_1 m^{-2}$ for some constant $K_1>0$ independent of $m$, the choice of $M$ and
the selection points  $\unu$, so that, for any $j$ and selection of $M$  and $\unu$, 
\be  \label{lower_bound}
\Delta_1 (j)  \geq K_2 2^{2j},
\ee
for some constant $K_2 >0$, independent of $j$. 
It turns out that if $M=M_n$ increases at least as fast as $n^{1/3}$, then, by sampling $u_l$, $l=1,2, \ldots, M$, 
uniformly on $U$, i.e., by selecting $u_l=a+(b-a)l/M$, $l=1,2,\ldots,M$,
one can attain $\Delta_1 (j)  \leq K_3 2^{2j}$ for some constant $K_3>0$, independent of $j$,
so that the upper and the lower bounds in this case coincide up to a constant independent of $n$
(see Pensky \& Sapatinas (2010)).

Unfortunately, the above results do not hold for finite values of $M$ or when $M=M_n$ is a slowly growing function of $n$. 
Indeed, in the case of small values of $M$, both $\tau_1 (m; \underline{u},M)$ and $\Delta_1 (j)$ have completely different dynamics from large $M$.
Indeed, if $M=1$, Johnstone \& Raimondo  (2004) and Johnstone {\em et al.}  (2004) showed that 
in the case of    $\gamma(u)=1/u$, $u_1=u^* = a = b$, one has $\Delta_1 (j)  \geq K_4 2^{3j}$ for any choice of   $u^*$ 
and some constant $K_4 >0$, independent of $j$. Johnstone {\em et al.}  (2004) also demonstrated that 
if $u^*$ is  selected to be a BA number, then the lower bound for $\Delta_1 (j)$ is attainable, i.e.,
$\Delta_1 (j)  \leq K_5 2^{3j}$ for some constant $K_5 >0$, independent of $j$. Hence, in this case, $\Delta_1(j) \asymp 2^{3 j}$.

These results were extended by De Canditiis \& Pensky (2006) who studied the multichannel deconvolution
model with a boxcar kernel and showed that the convergence rates obtained by Johnstone {\em et al.}  
(2004) for $M=1$ can be improved by sampling at several different points. In particular, they demonstrated that 
if $M$
is \textit{finite}, $M \geq 2$, one of the $u_1,u_2,\ldots,u_M$ is a BA
number, and $\underline{u}$ is a BA 
$M$-tuple defined in \fr{BAtuple}, then  
\be \label{infMbox}
\Delta_1 (j) \leq C(M) \, j 2^{j(2+ 1/M)}
\ee
for some positive  $C(M)$. 
In particular, when $M$ is growing with $n$, the value of $C(M)$ 
depends on $n$ and, hence, affects the convergence rates of the estimator $\hf(\cdot)$ as $n \rightarrow \infty$.

The relation between the convergence rates of the estimator $\hat{f}_n(\cdot)$, given by \fr{fest}, of $f(\cdot)$
in the model \fr{eq:box-car1} and the theory of Diophantine approximation
becomes obvious when one notes that in \fr{eq:box-car-taulower}, for any $m \in \ints \setminus \{0\}$ and any $u_l$, $l=1,2,\ldots,M$, one has, combining the periodic behavior of the sine function together with a first order (linear) approximation,
\be  \label{sin_dioph}
4 \|2 m u_l\|^2 \leq \sin^2 (2 \pi m u_l)   \leq \pi^2\, \|2 m u_l\|^2,\ \ l=1,2,\ldots, M,
\ee
where $\| a \|=\inf\{|a-k|,\, k \in \ints\}$ denotes  the distance from a real number $a$ to the nearest integer number. Hence, \fr{eq:box-car-taulower} becomes
\begin{align}
\label{eq:box-car-diophant}
\tau_1 (m) \equiv \tau_1 (m; \underline{u},M) & 
\asymp \frac{1}{m^2\, M}   \sum_{l=1}^M \|2 m u_l\|^2,
\end{align}
so that the convergence rates of the estimator $\hat{f}_n(\cdot)$ depend on the lower bound, in terms of $m$, of the expression
\fr{eq:box-car-diophant}.

The value of $C(M)$ in \fr{infMbox} is related to the value of $B(M)$ in \fr{BAtuple}.
To the best of our knowledge, there has not been developed a procedure for construction
of a BA $M$-tuple on a specified interval, of non-asymptotic length, of the real line,  and there are no asymptotic lower bounds, in terms of $M$, 
on $B(M)$ in \fr{BAtuple} when the value of $M$ is growing. 
For this reason, in order to find upper bounds of estimator \fr{fest} and choose an optimal  relation between
the sample size $n$ and the number of channels $M$ when $M=M_n$ is a slowly growing function of $n$,
we need to obtain new original results in Diophantine approximations. In particular, 
the objective of the next section  is to construct a BA $M$-tuple on the non-asymptotic interval $U$, of the real line,
and to obtain a lower bound on $B(M)$ in terms of $M$ for this BA M-tuple when $M$ grows slowly with $n$.

\section {Construction of a BA ${\bf M}$-tuple on a specified interval }
\label{sec:num_theory}
\setcounter{equation}{0}

Below, we construct   a BA  $M$-tuple $\ubeta= (\beta_1, \beta_2, \ldots,  \beta_M)$ of real numbers on a specified interval $(a,b)$,
of a non-asymptotic length, of the real line, and derive the  lower bound on $B(M)$ in  formula \fr{BAtuple}. 
For this construction, we use the technique described in Schmidt (1980), Section 4, pp. 43-45.
In particular, we shall provide an algorithm for construction of an $M$-tuple
$\beta_1, \beta_2, \ldots,  \beta_M$ of real numbers such that, as $M \rightarrow \infty$, 
\begin{enumerate}
\item it lies in any specified interval $(a,b)$, $a < b$, of nonasymptotic length, of the real line, and
\item it satisfies
\be  \label{BAtuple Lower}
\max_{i=1,2,\ldots, M} | \beta_i q - p_i | \geq B_0 \exp (-6M \ln M) q^{-1/M}, 
\ee
for any integer numbers  $q>0$ and $p_1, p_2, \ldots, p_M$,
and for some constant $B_0>0$, independent of $M$, 
$q$ and $p_1, p_2, \ldots, p_M$, so that $B(M) = B_0 \exp (-6M \ln M)$ in \fr{BAtuple}.
\end{enumerate}

Assume that $M$ is large enough, fix a positive integer $Q$ and consider 
\be  \label{polynom}
P(x) = (x-Q)(x-2Q)\cdots (x-MQ)-1,
\ee
a monic polynomial (i.e., a polynomial with a unit leading coefficient) of the degree $M$. 
Let $\xi_1, \xi_2, \ldots, \xi_M$ be the roots of a polynomial \fr{polynom}.
Recall that $\xi$ is called an {\em algebraic integer} number if it is a root of some monic polynomial   with coefficients being integer numbers.
Algebraic integers are called {\it conjugate} if they are roots of the same monic polynomial with integer coefficients.

\medskip

Then, the following statement is valid.

\begin{lemma}  \label{lem2}
If  $Q \geq 5M$, then $\xi_1, \xi_2, \ldots, \xi_M$ are real conjugate algebraic integer numbers 
such that
\be \label{alg_int}
(i-1/2)Q < \xi_i < (i+1/2)Q,\ \  i=1,2, \ldots, M.
\ee
\end{lemma}

Now, to construct a BA $M$-tuple, choose $Q \geq 5(M+1)$ and construct real conjugate algebraic integers  $\xi_1, \xi_2, \ldots, \xi_M, \xi_{M+1}$ 
using the process described in Lemma \ref{lem2}.
Let $\uaf = (\af_1, \af_2, \ldots, \af_M)$ be a solution to the following system of equations:
\be \label{largeAI}
\sum_{i=1}^M \xi_k^{i-1} \af_i = - \xi_k^M,\  \  k=1,2, \ldots, M.
\ee 
Observe that the determinant of the system of equations \fr{largeAI} 
is  a Vandermonde determinant; hence, it is nonzero since  $\xi_i \neq \xi_j$ for $i \neq j$.
Therefore, the system of linear equations \fr{largeAI} has a unique  solution $\uaf = (\af_1, \af_2, \ldots, \af_M)$,
which turns out to be a BA $M$-tuple.
 
\begin{lemma}  \label{lem3}
The solution $\uaf = (\af_1, \af_2, \ldots, \af_M)$ of the system of equations \fr{largeAI} is a 
BA $M$-tuple such that
\be \label{bounds}
|\af_k| \leq 30 \exp(3M \ln M), \ \ k=1,2, \ldots, M, 
\ee
and for any integer numbers $q>0$ and $p_1, p_2, \ldots, p_M$, as $M \rightarrow \infty$,  one has 
\be  \label{BAbounds}
\max_{i=1,2, \ldots, M} | \af_i q - p_i | \geq C_0 \exp (- 3 M \ln M)\,  q^{-1/M}  
\ee
with some constant $C_0 >0$, independent of $M$, $q$ and $p_1, p_2, \ldots, p_M$.
\end{lemma}

Lemma \ref{lem3} provides a BA $M$-tuple which, however, does not necessarily belong to the specified interval $(a,b)$,
of a non-asymptotic length, of the real line.
Assume, without loss of generality, that both $a$ and $b$ are rational numbers, otherwise, replace $(a,b)$
by $(a^*, b^*) \in (a,b)$, where $a^*$ and $b^*$ are rational numbers. Let $a=p_a/q_0$ and $b=p_b/q_0$ for  
some integer numbers  $p_a$, $p_b$ and $q_0$,
and let $z$ be an integer number such that $z-1 < 30\, \exp(3M \ln M) \leq z$. 
Define
\be \label{smallBA}
\beta_l = a +  \af_l (b-a)/z, \quad l=1,2,\ldots,M,
\ee
where $\uaf = (\af_1, \af_2, \ldots, \af_M)$  is the BA $M$-tuple constructed in Lemma \ref{lem3}.

\medskip

The following theorem confirms that $\ubeta= (\beta_1, \beta_2, \ldots,  \beta_M)$, as constructed above, forms a BA $M$-tuple
on the specified $(a,b)$, of a non-asymptotic length, of the real line.
 
\begin{theorem}  \label{th1}
The real numbers $\beta_1, \beta_2, \ldots, \beta_M$ defined in \fr{smallBA} lie on the interval
$(a,b)$, of a non-asymptotic length,  and form a BA $M$-tuple, so that, as $M \rightarrow \infty$, one has
\be  \label{FinalBA}
\max_{i=1,2, \ldots, M} | \beta_i q - p_i | \geq B_0 \exp (- 6 M \ln M)\,  q^{-1/M},  
\ee
for any integer numbers 
$q>0$ and $p_1, p_2, \ldots, p_M$,  and for some constant $B_0 >0$, independent of $M$, $q$ and
$p_1, p_2, \ldots, p_M$, so that $B(M) =  B_0 \exp (- 6 M \ln M)$.
\end{theorem}

\section {Asymptotical  upper bounds for the $L^2$-risk of  the adaptive wavelet thresholding estimator }
\label{sec:asymp_bounds}
\setcounter{equation}{0}

In Section \ref{sec:num_theory}, we constructed a BA $M$-tuple and derived a lower bound on $B(M)$ in \fr{BAtuple}, as $M \rightarrow \infty$.
We can now choose the resolution levels $j_0$ and $J$, the thresholds $\lam_{j}$ in \fr{fest} and the optimal relation between the 
total number of observations $n$ and the number of channels $M=M_n$ and derive asymptotical upper bounds for the $L^2$-risk  
of the estimator $\hf(\cdot)$ given by \fr{fest} over a collection of Besov balls.

In order to formulate and prove Theorem \ref{th:box-car-slow} we first need to obtain some preliminary results. 
Recall that  $\| a \|$ denotes  the distance from a real number $a$ to the nearest integer number.
For this purpose, we recall the   {\em equidistribution} lemma (see Lemma \ref{lem1}), 
proved in Johnstone  \&  Raimondo  (2004), which we state here for completeness, and formulate a new lemma
(i.e., Lemma \ref{lem_aleph}) which is based on application of Lemma \ref{lem1} to the BA $M$-tuple.

\begin{lemma} \label{lem1} (Lemma 1 in Johnstone  \&  Raimondo  (2004)) 
Let $p/q$ and $p'/q'$ be successive convergents in the continued fraction expansion of a real number $a$. Let $N$ be a positive integer number with $N + q < q'$.
Let $h$ be a non-increasing function. Then 
$$ 
\sum_{i=4}^q h(i/q) \leq \sum_{k=N+1}^{N+q} h(\| ka \|) \leq 2 \sum_{i=1}^{q-3} h(i/q)
+ 6 h(1/(2q')).
$$
\end{lemma}

\begin{lemma}  \label{lem_aleph}
Let $\beta_1, \beta_2, \ldots, \beta_M$ be a BA $M$-tuple constructed in Theorem \ref{th1} and
let $\beta_1$ be a BA number. Let $r_0$ be an arbitrary fixed positive real number. Denote
\be  \label{aleph}
\aleph_k (j,M)  = \sum_{l \in \Om_j} \lkv \| l \beta_1 \|^2 + \cdots +  \| l \beta_M \|^2 \rkv^{-k},
\ee
where  $\Om_j$ is defined as
\be  \label{Omj}
\Om_j = \lfi l: 2^j \leq |l| \leq 2^{j + r_0} \rfi.
\ee
If $M$ is large enough, then, as $j \rightarrow \infty$, 
\be \label{aleph_bounds}
\aleph_k (j,M) = O\lkr j\ 2^{j(1 + (2k-1)/M)}\ e^{6(2k-1) M \ln M} \rkr, \quad k=1,2,3,4.
\ee
\end{lemma}

We also need the following two lemmas which evaluate the precision of estimation of $\ajk$ and $\bjk$.

\begin{lemma}  \label{lemm:bound1}
Let $\ubeta= (\beta_1, \beta_2, \ldots, \beta_M)$ be a BA $M$-tuple constructed on the interval $(2a, 2b)$ 
according to Theorem \ref{th1} and let one of $\beta_1, \beta_2, \ldots, \beta_M$ be a BA number.  
Let the equation \fr{eq:box-car1} be evaluated at the 
the point $\underline{u}$ with components $u_l=\beta_l/2$, $l=1,2,\ldots,M$. 
Then,  for all $j \geq j_0$, as $n \rightarrow \infty$,
\begin{align*}
\EE |\hajk - \ajk|^2   &  = O \lkr n^{-1} {j_0}\ M \ 2^{{j_0}(2 + 1/M)}\   e^{6M \ln M}  \rkr, \\
\EE |\hbjk - \bjk|^{2} &  = O \lkr n^{-1} j\ M \ 2^{j(2 + 1/M)}\   e^{6M \ln M} \rkr, \\
\EE |\hbjk - \bjk|^{4} &  = O \lkr n^{-2} \ j^2\ M^2\ 2^{j(4 + 2/M)}\ e^{42 M \ln M}  \rkr.
\end{align*}
\end{lemma}

\begin{lemma}  \label{lemm:bound2}
Let $u_1, u_2, \ldots, u_M$ be as  in Lemma \ref{lemm:bound1}.  If $\eta >0$ is 
a constant large enough, then, for all $j \geq j_0$, as $n \rightarrow \infty$,
$$
\PP (|\hbjk - \bjk|^2 \geq \eta^2 (n_M)^{-1} j 2^{j(2+1/M)} \ln n)=o \lkr n^{-\theta}\rkr,
$$
where $n_M=(n/M) \exp(-6M\ln M)$ and $\theta=\eta^2/(2C_{\psi})$ with $C_{\psi}=2^{-j}|C_j|$.
\end{lemma}

We are now ready to formulate Theorem \ref{th:box-car-slow}.
Let the resolution levels $j_0$
and $J$ and the thresholds $\lam_{j}$ be such that
\be \label{ff:newparam}
2^\jo = \ln n, \ \ 2^J = (n_M)^{\frac{1}{3+1/M}}, \ \
\lam_j = \eta \; (n_M)^{-1/2}\, \sqrt{j 2^{j(2+1/M)}\ln n},
\ee
for some constant $\eta >0$, where
\be
\label{eq:nM}
n_M = \frac{n}{M} \exp(-6 M \ln M).
\ee
Note that since the construction of $j_0$, $J$ and $\lambda_j$ is independent of the Besov ball parameters, $s$,
$p$, $q$ and $A$, the suggested wavelet thresholding estimator $\hf(\cdot)$ given by \fr{fest}
is adaptive with respect to these parameters.
\\

 The  following statement provides the asymptotical upper bounds for the
$L^2$-risk, over a collection of Besov balls.

\begin{theorem} 
\label{th:box-car-slow}
Let $s >1/\min(p,2)$, $1 \leq p \leq \infty$, $1 \leq q \leq \infty$ and
$A>0$. Let $\ubeta= (\beta_1, \beta_2, \ldots, \beta_M)$ be a BA $M$-tuple constructed on the interval $(2a, 2b)$ 
according to Theorem \ref{th1} and let one of $\beta_1, \beta_2, \ldots, \beta_M$, say $\beta_1$, be a BA number. 
Let the equation \fr{eq:box-car1} be evaluated at the 
the point $\underline{u}$ with components $u_l=\beta_l/2$, $l=1,2,\ldots,M$. Choose 
\be  \label{Mnval}
M = M_n =\nu\, \sqrt{\ln n/(\ln \ln n)}
\ee 
for some $\nu \leq 1/\sqrt{6}$, independent of $n$. 
Let $\hat{f}_n (\cdot)$ be the adaptive wavelet thresholding estimator
defined by \fr{fest} with $\jo$, $J$ and $\lambda_j$  given by \fr{ff:newparam}, and $n_M$ given by \fr{eq:nM}. 
Then,  as $n \rightarrow \infty$,
\begin{equation}
\label{eq:box-car-lowerF1SLOWLY} 
R_n(\hat{f}_n)  \leq \lfi
\begin{array}{ll}
C \, n^{-\frac{2s}{2s+3}}\,  \alpha_n, & \mbox{if}\;\;\; s
>3(1/p - 1/2),
\\
C \, \Big ( \frac{\ln n}{n} \Big )^{\frac{s'}{s'+1}} \,  \alpha_n, & \mbox{if}\;\;\; s \leq 3(1/p - 1/2),
\end{array} \right.
\end{equation}
where  $\alpha_n$ is given by
\begin{equation}
\label{Newalphan}
\alpha_n = \exp \lfi \sqrt{\ln n} \, \sqrt{\ln \ln n} \lkv \frac{A_1}{A_2} \lkr 3\nu + \frac{1}{A_2 \nu} \rkr + r_n \rkv\rfi,
\end{equation}
with 
\begin{eqnarray*}
r_n & = & 
\frac{3 A_1\nu\ \ln\ln\ln n}{A_2\ \ln \ln n}\lkr \frac{2\, \ln \nu}{\ln\ln\ln n} -1 \rkr 
+ \frac{\sqrt{\ln \ln n}}{\sqrt{\ln n}} \lkr \frac{A_3}{A_2} +  \frac{A_1}{2 A_2} - \frac{3 A_1}{A_2^2} - \frac{A_1}{A_2^3 \nu^2} \rkr\\  
& + & 
\frac{\ln \ln \ln n}{\sqrt{\ln n} \sqrt{\ln \ln n}} \lkr \frac{3 A_1}{A_2^2} - \frac{A_1}{2 A_2} \rkr = o(1) \,\, (n \rightarrow \infty),
\end{eqnarray*}
where 
\begin{equation}  \label{A123}
\begin{array}{ll}
A_1=2s, \,\,\,\,A_2 = 2s+3, \,\,\,\,\,A_3=2s, & \mbox{if}\;\;\;  2 \leq p \leq \infty,  \\
A_1=2s, \,\,\,\,A_2 = 2s+3, \,\,\,\,\,A_3=4s, & \mbox{if}\;\;\; 6/(2s+3) < p < 2,   \\
A_1=2s^*, \,A_2 = 2s^*+3, \,\,A_3=4s^*, & \mbox{if}\;\;\; 1 \leq p \leq 6/(2s+3)   
\end{array}
\end{equation}

with $s^*=\min(s',s)$, $s'=s+1/2-1/p$.
\end{theorem}

\section {Discussion }
\label{sec:discussion}
\setcounter{equation}{0}

We considered the estimation problem of the unknown response function in the multichannel boxcar deconvolution model with a boxcar-like kernel when the number of channels grows 
as the total number of observations increases.  This situation seems to be of a particular interest nowadays  since data recording equipment is getting cheaper and cheaper while overall volumes of data is growing very fast. Our aim was to investigate the situation when the number of channels $M=M_n$ grows slowly with the number of observations $n$. 

For this purpose, we obtained new original results  in the field of Diophantine approximation in order to devise a technique which allows the reconstruction of the unknown response function with a precision that differs from the best possible convergence rates (which can be attained in the corresponding continuous functional deconvolution model \fr{convcont}) by a factor which grows slower than any power of $n$.

Specifically, in Section \ref{sec:asymp_bounds}, we derived asymptotical upper bounds for the $L^2$-risk of the adaptive wavelet thresholding estimator
\fr{fest} of $f(\cdot) \in L^2(T)$ in the model \fr{eq:box-car1}. In comparison, it  follows from Pensky \& Sapatinas (2010) that
the choice of a uniform sampling strategy (i.e., $u_l=a+(b-a)l/M$, $l=1,2,\ldots,M$, for $M =M_n \geq  (32 \pi/3)(b-a)n^{1/3}$,  
leads to an adaptive wavelet block thresholding estimator $\hat{f}_n^{\rm B}(\cdot)$ of $f(\cdot)$ with the following convergence rates 
\begin{equation} \label{box-car-fast} 
R_n (\hat{f}_n^{\rm B})  \leq \lfi
\begin{array}{ll}
C n^{-\frac{2s}{2s+3}}\ \lkr  \ln n\rkr^{\ro}, & \mbox{if}\;\;\; s>3(1/p - 1/2),
\\
C \Big ( \frac{\ln n}{n} \Big )^{\frac{s'}{s'+1}} \ \lkr  \ln n
\rkr^{\ro}, & \mbox{if}\;\;\; s \leq 3(1/p - 1/2),
\end{array} \right.
\end{equation}
for $s > 1/\min(p,2)$, where $s' = s+1/2-1/p$ and 
$$ \ro = \lfi
\begin{array}{ll}
\frac{3\max(0, 2/p-1)}{2s+3}, & \mbox{if}\;\;\; s
>3(1/p - 1/2), \\  
\max(0, 1 -p/q), & \mbox{if}\;\;\; s = 3(1/p - 1/2),  \\
0, & \mbox{if}\;\;\; s < 3(1/p - 1/2).   
\end{array} \right.
$$
Moreover, its has been shown that the above convergence rates with $\rho=0$ are the fastest possible ones (see Pensky \& Sapatinas (2010)); 
hence,  up to the logarithmic factor $\lkr \ln n\rkr^{\ro}$, $\hat{f}_n^{\rm B}(\cdot)$ attains the best possible convergence rates.
By comparing the convergence rates \fr{eq:box-car-lowerF1SLOWLY} in Theorem \ref{th:box-car-slow} to the fastest possible convergence rates 
(without the extra logarithmic factor $\lkr \ln n\rkr^{\rho}$ appearing in (\ref{box-car-fast})), one concludes that they differ by the extra factor $\alpha_n$ defined in \fr{Newalphan}.

How fast does $\alpha_n \rightarrow \infty$ as $n \rightarrow \infty$? 
It can be easily seen that  $\alpha_n$ grows slower than any power of $n$ but faster than any power of $\ln n$, i.e.,
for any $a_1, a_2 >0$,   one has
$$
\lim_{n \rightarrow \infty} \frac{\alpha_n}{n^{a_1}}=0, \quad \lim_{n \rightarrow \infty} \frac{\alpha_n}{(\ln n)^{a_2}}= \infty.
$$
Hence, although choosing $M=M_n \rightarrow \infty$ at a rate given by \fr{Mnval} improves the convergence rates in comparison
with the finite values of $M$ (see Pensky \& De Canditiis (2006), Theorem 2), these rates are quite a bit worse than in the case when $M=M_n$ grows at a faster rate as $n \rightarrow \infty$.
Since, as we have explained in Section \ref{sec:relation}, this fast growth of $M=M_n$ with $n$ cannot be achieved in a number 
of practical situations, one has to resign  to $M_n$ growing slowly with $n$, in particular, 
$M=M_n =o((\ln n)^{\alpha_3})$ for some $\alpha_3 \geq 1/2$.

The interesting question, however, is whether the  convergence rates \fr{eq:box-car-lowerF1SLOWLY}  can be improved.
To uncover an answer to this question, one needs to either come up with another procedure for constructing a BA $M$-tuple  which belongs 
to a  specified interval, of a non-asymptotic length, of the real line and delivers a higher value of $B(M)$ in \fr{BAtuple}, as $M \rightarrow \infty$, or to show 
that no matter what the value of $\unu = (u_1, u_2, \ldots, u_M)$ is, there exist integer numbers $q$ and $p_1,p_2,\ldots,p_M$ such that, as $M \rightarrow \infty$, 
$$
\max_{i=1,2,\ldots, M} | u_i q - p_i | \leq B_1 \exp (-6M \ln M) q^{-1/M}, 
$$
for some positive constant $B_1$ independent of $M$, $q$ and $p_1,p_2,\ldots,p_M$. At the moment we are unable to provide answers to either of the above questions; we challenge, however,
the number theory community to work on the issue. Derivation of these results will not only enrich the theory of Diophantine approximation
but will also be valuable for the theory of statistical signal processing.

\section {Proofs }
\label{sec:proofs}
\setcounter{equation}{0}

 \noindent {\bf Proof of Lemma \ref{lem2}.} Observe that for $P(x)$ given by \fr{polynom}, one has
$$
P((M+1/2)Q)>0,\ \  (-1)P((M-1/2)Q)>0,\ \  \ldots, \ \ (-1)^M P(Q/2)>0,
$$
so that $P(x)$ has $M$ real roots $\xi_1, \xi_2, \ldots, \xi_M$ such that 
\fr{alg_int} is valid. By definition, $\xi_1, \xi_2, \ldots, \xi_M$ are algebraic integer numbers.
Let us show that no proper subset of $\xi_1, \xi_2, \ldots, \xi_M$ is itself a set
of conjugate algebraic integer numbers. For this purpose, note that 
\be \label{farroot}
Q(|j-i|  - 1/2) \leq |\xi_i -jQ| \leq Q(|j-i|  + 1/2), \ \ i \neq j,\ i,j=1,2, \ldots, M. 
\ee
Therefore, by \fr{polynom}, for any $i=1,2, \ldots, M,$
\be  \label{closeroot}
0 \leq |\xi_i - iQ| = \lkv \prod_{\stackrel {j=1}{j\neq i}}^M |\xi_i - jQ| \rkv^{-1} 
\leq Q^{-(M-1)} \prod_{\stackrel {j=1}{j\neq i}}^M (|j-i|-1/2)^{-1}.
\ee
Now, assume that $\xi_{i_1}, \xi_{i_2}, \ldots, \xi_{i_m}$, $i_1 < \ldots < i_m$ and $m<M$, form a set of conjugate real integer numbers. 
Then, $P^*_m (i_1 Q) = (\xi_{i_1} - i_1 Q) \ldots  (\xi_{i_m} - i_1 Q)$ is an integer number and is not equal to zero, hence, 
$|P^*_m (i_1 Q)| \geq 1$. On the other hand, by \fr{farroot} and \fr{closeroot}, 
$$
1 \leq |P^*_m (i_1 Q)| \leq Q^{-(M-1)} \prod_{\stackrel {j=1}{j\neq i_1}}^M (|j-i_1| -1/2)^{-1}\ \prod_{k=2}^m [Q(|i_k - i_1|+1/2)].
$$
The product in the right-hand side above takes the largest value if $m = M-1$, $i_1=1$ and $i_k=k+1$, $k=2,3, \ldots, M-1$.
In this case, for $M>2$, combination of the last two inequalities yields
$$
1 \leq |P^*_m (i_1 Q)| \leq 4 (M-1/2) Q^{-1} < 5M/Q,
$$
which leads to a contradiction when $Q >5M$.  \hfill $\Box$

\bigskip

\noindent {\bf Proof of Lemma \ref{lem3}.} Choose $Q = 5(M+1)$ and construct real conjugate algebraic integer numbers  $\xi_1, \xi_2, \ldots, \xi_M, \xi_{M+1}$ 
using the process described in Lemma \ref{lem2}. Then, by \fr{closeroot},  $\xi_i \approx Q i$.
Let $q>0$ and $\unp =(p_1, p_2, \ldots, p_M)$ be integer numbers and denote
$$
H_k (q, \unp) = \sum_{i=1}^M \xi_k^{i-1}  p_i + \xi_k^M q, \ \ k=1,2, \ldots, M+1.
$$
Note that if $\unp$ is not zero and the components of the vector $\unp/q$ are not integer
numbers, then $H_k (q, \unp) \neq 0$. 
Furthermore, $H_1 (q, \unp), H_2 (q, \unp), \ldots, H_M (q, \unp), H_{M+1} (q, \unp)$ are themselves real conjugate
algebraic integer numbers  and, thus, 
$$
\prod_{k=1}^{M+1} |H_k (q, \unp)| \geq 1.
$$
Now, note that \fr{largeAI} implies that $\af_1, \af_2,\ldots, \af_M$ are coefficients of the monic polynomial with the roots 
$\xi_1, \xi_2, \ldots, \xi_M$. Also, it is easy to check that for the solution  $\uaf=(\af_1, \af_2,\ldots, \af_M)$ of the system of 
equations \fr{largeAI}, one can write
\be  \label{hkq}
H_k (q, \unp) = \sum_{i=1}^M \xi_k^{i-1} (p_i - \af_i q), \ \ k=1,2, \ldots, M.
\ee
Moreover, if one denotes
\be \label{omegaM}
\om_M = \sum_{i=1}^{M} \af_i \xi_{M+1}^{i-1} + \xi_{M+1}^M,
\ee
then $H_{M+1} (q, \unp)$ can be written as
\be  \label{hmqp}
H_{M+1} (q, \unp) = \sum_{i=1}^M \xi_{M+1}^{i-1} (p_i - \af_i q) + \om_M.
\ee
Note that \fr{omegaM} implies that $\om_M$ is the value of the polynomial 
$$
{\cal P}(x) = \sum_{i=1}^M \af_i x^{i-1} + x^M = (x-\xi_1) \cdots (x-\xi_M)
$$ 
at the point $\xi_{M+1}$. Therefore, by \fr{farroot} and \fr{closeroot}, 
$|\om_M| \leq K M! Q^M$ for some constant $K>0$. 

Recall that $\| a \|$ denotes the distance from a real number $a$
to the nearest integer number. Denote $L= \max_{i=1,2,\ldots,M}|\xi_i q - p_i|$.
Note that we can assume that $p_i/q$, $i=1,2, \ldots, M$, are not integer numbers. Otherwise,
if, for instance, $p_1/q=z$ is an integer number, then
$L \geq q |\xi_1 - z| \geq q \| \xi_1 \|$  and \fr{BAbounds} is valid. 
If $L \geq 1$, then \fr{BAbounds} is valid. Hence, consider the case of $L<1$. Then $L <|q|$
and, by \fr{hkq}, we have 
$$  
|H_k (q, \unp)| \leq L \sum_{i=1}^M \xi_k^{i-1} < L \xi_k^M/(\xi_k-1), \ \ \ k=1,2, \ldots, M.
$$
Then, using \fr{hmqp} and an upper bound for $\om_M$, we obtain
$$
|H_{M+1} (q, \unp)| \leq |q| (\xi_{M+1}^M/(\xi_{M+1}-1) + K M! Q^M).
$$
Since $H_1 (q, \unp), H_2 (q, \unp), \ldots, H_M (q, \unp), H_{M+1} (q, \unp)$ are  real conjugate
algebraic integer numbers, one has
\begin{equation*}  \label{invertform}
1 \leq \prod_{k=1}^{M+1} |H_k (q, \unp)| \leq \prod_{k=1}^{M} \lkv \frac{L \xi_k^M}{\xi_k-1} \rkv |q| 
\lkv \frac{\xi_{M+1}^M}{\xi_{M+1}-1} + K M! Q^M \rkv.
\end{equation*}
Note that, by \fr{closeroot} and  $Q \geq 5M$,  one has $|\xi_{M+1} - (M+1)Q | <Q^{-M}$
and, hence, $|\xi_{M+1}|^M \leq 2 Q^M (M+1)^M$. Therefore
$$
L \geq K |q|^{-1/M} \ \prod_{k=1}^M \lkv \frac{Qk-1}{k^M Q^M}  \rkv^{1/M} \ 
\lkv \frac{(M+1)^M Q^M}{Q(M+1)-1} + Q^M M! \rkv^{-1/M}.
$$
Plugging in $Q=5(M+1)$ into the expression above, we obtain 
$$L \geq B(M) |q|^{-1/M}$$ 
with 
$$
B(M) = K \, [5(M+1)]^{-(M+1)} (M!)^{-1}  \prod_{k=1}^M [5k(M+1) -1]^{1/M}\ 
\lkv \frac{(M+1)^M }{5(M+1)^2 -1} +   M! \rkv^{-1/M}.
$$
Using Stirling formula,
$$M! =  \sqrt{2 \pi} (M+1)^{M +1/2} \exp(-(M+1)) (1 + o(1)), \quad \text{as} \quad M \rightarrow \infty,$$
(see, e.g., formula 8.327  of Gradshtein and Ryzhik  (1980)) and the fact that $\ln (M+1) < \ln(M) + 1/M$, after some simple algebra,
we obtain that, as $M \rightarrow \infty$, $$B(M) \geq C_0 \exp(- 3M \ln M),$$ 
for some constant $C_0 >0$, independent of $M$, $q$ and $\unp$, which proves \fr{BAbounds}.

Now, it remains to prove  the upper bound  \fr{bounds} for $\af_k$, $k=1,2,\ldots,M$.
For this purpose, recall that $\af_1, \af_2, \ldots, \af_M$ are coefficients of the monic polynomial with 
roots $\xi_1, \xi_2, \ldots, \xi_M$. Therefore, using  \fr{closeroot}, obtain
$$
|\af_k| \leq  {M \choose k} M Q (M-1)Q \cdots (M-k+1)Q = k! {M \choose k}^2 Q^k, \quad k=1,2,\ldots,M.
$$
Since for any $k=1,2,\ldots,M$, $5^k/k! \leq 625/24 <30$,  $Q=5(M+1)$ and $(M+1)(M-j) \leq M^2$ for $j \geq 1$, one has  
(reading $\prod_{j=0}^{-1} =1$)
\beqns
|\af_k| & \leq  &  \frac{5^k}{k!}  \ (M+1)^k \prod_{j=0}^{k-1}  (M-j)^2 \\
        & \leq & 30\ M^2 [(M+1)(M-k+1)]^2 \ \prod_{j=1}^{k-2} [(M+1) (M-j)^2] \\
        & \leq & 30 M^{3k} \leq 30 e^{3M \ln M}, \quad k=1,2,\ldots,M,
\eeqns
which proves \fr{bounds}. \hfill $\Box$

\bigskip

\noindent {\bf Proof of Theorem \ref{th1}.} It is easy to check that $\beta_1,\beta_2,\ldots,\beta_M$, as defined
by \fr{smallBA}, lie on  $(a,b)$. Furthermore,  by Lemma \ref{lem3} 
and the fact that  $z < 30\,  \exp(3M \ln M)$, as $M \rightarrow \infty$, one has
\beqns
\max_{i=1,2,\ldots, M} | \beta_i q - p_i | & = & 
(z q_0)^{-1} \,  \max_{i=1,\cdots, M} |\af_k (p_b - p_a)q - (z q_0 p_l - z p_a q)|\\
& \geq &  (z q_0)^{-1} \, C_0 |(p_b - p_a)q|^{-1/M} \exp(- 3M \ln M) \\
& \geq&  B_0 \exp(- 6M \ln M)\ |q|^{-1/M}, 
\eeqns
for any integer numbers 
$q>0$ and $p_1, p_2, \ldots, p_M$,  and for some constant $B_0 >0$, independent of $M$, $q$ and
$\unp$. \hfill $\Box$

\bigskip

\noindent {\bf Proof of Lemma  \ref{lem_aleph}.}  
Recall first that any real number $a$, which is not an integer number, may be uniquely determined by its
continued fraction expansion
$$
a = a_0 + \frac{1}{a_1 + \frac{1}{a_2 + \frac{1}{a_3 + \cdots}}},
$$
where $a_0$ is an integer number and $a_1, a_2, \ldots$ are strictly
positive integer numbers. The convergents $p_k/q_k=p_k(a)/q_k(a)$, $k=0,1,\ldots$,
of $a$ are those rational numbers, the continued fraction expansions of
which terminate at stage $k$, that is, $p_0/q_0 = a_0$, $p_1/q_1
=a_0 + 1/a_1$, $p_2/q_2 = a_0 + 1/(a_1 + 1/a_2)$, and so on. The
denominators in the above expansions grow at least geometrically
\begin{eqnarray}
q_{n+i}  &\geq & 2^{(i-1)/2} q_n, \quad \mbox{if}\;\;\;  i\ \mbox{odd},  \label{dioph1}\\
q_{n+i}  & \geq & 2^{i/2} q_n,  \quad \quad \;\; \mbox{if}\;\;\;  i\ \mbox{even},\nonumber 
\end{eqnarray} 
and  $a_n < q_n/q_{n-1} \leq a_n + 1$, $n \geq 1$.
A real number $a$ is BA if $\sup_n   a_n  < \infty$, i.e.,
there exists $\barq >0$ such that
\be
q_n/q_{n-1} \leq \barq,\ \ n \geq 1
\label{dioph2}
\ee
(see, e.g., Schmidt (1980), Sections 3-5, pp. 7-23).

Let $p/q$ and $p'/q'$ be successive principal
convergents in the continued fraction expansion of 
$\beta_1$.  Let $N$ be a positive integer number with $N + q < q'$.
Then, application of Lemma \ref{lem1} with $h(x) = x^{-1}$ 
yields
\be  \label{powerrel}
\sum_{l=N+1}^{N+q} \| l \beta_1 \|^{-1} = O( q \ln q ),\ \ \ 
\ee
since $q' \leq \barq q$ by \fr{dioph2}. Now, note that by \fr{FinalBA}
\be  \label{sumequa}
\sum_{l=N+1}^{N+q} \lkr \| l \beta_1 \|^2 + \cdots +  \| l \beta_M \|^2 \rkr^{-k} \leq 
\sum_{l=N+1}^{N+q} \| l \beta_1 \|^{-1} \ [\max(\| l \beta_1 \|, \ldots , \| l \beta_M \|)]^{-(2k-1)}\  
\ee
Combination of \fr{FinalBA}, \fr{powerrel}  and \fr{sumequa} implies that 
\be \label{inequlity1} 
\sum_{l=N+1}^{N+q} \lkr \| l \beta_1 \|^2 + \cdots +  \| l \beta_M \|^2 \rkr^{-k} 
=  O  \lkr e^{6(2k-1)\, M \ln M}\ q^{(1 + (2k-1)/M)} \ln q \rkr,\ k=1,2,3,4.\ \ \ 
\ee

Now, observe that the set of indices $l$ in $\Om_j$ is symmetric
about zero, and so are the components of the sum. Hence, we can consider
only the positive part of $\Om_j$ which, with some abuse of notation, we keep calling it $\Om_j$.
Let $q_i$ be the denominators of the convergents of $\beta_1$, and let
$l$ be the smallest number such that $q_l \geq 2^j$. The geometric
grows of denominators \fr{dioph1} implies that $2^{j+r_0} <2^{r_0} q_l \leq q_{l+ 2r_0}$ so that $\Om_j  \subseteq [q_{l-1}, q_{l+ 2 r_0})$.
If we denote $D_s = N \cap [q_{l+s-1}, q_{l+s})$, $s=0,1, \ldots, 2r_0$, then
$$
\Om_j \subseteq \bigcup_{s=0}^{2 r_0} D_s.
$$
Since, by \fr{dioph1}, $q_{i+1} \leq \barq q_i$, there are at most $\barq$
disjoint blocks of length $q_{l+s-1}$ that cover $D_s$. 
Applying  \fr{inequlity1} to each of those blocks, we derive
$$
\sum_{l \in D_s} \lkr \sum_{i=1}^M  \|l \beta_i\|^2 \rkr^{-k} =
O \lkr  e^{6(2k-1)\, M \ln M}\, (q_{l+s-1})^{1 + (2k-1)/M}\  \ln q_{l+s-1}  \rkr,\ \ k=1,2,3,4.
$$
Note that $q_{l-1} \leq  2^j$, so that $q_{l+s-1} \leq \barq^s q_{l-1} \leq  \barq^s 2^j$.
Therefore,
\beqns
\aleph_k (j, M) & = & O \lkr   \sum_{s=0}^{2 r_0} \sum_{l \in D_s} \lkr \sum_{i=1}^M  \|l \beta_i\|^2 \rkr^{-k} \rkr
\nonumber \\
& = & O \lkr  e^{6(2k-1)\, M \ln M}  \sum_{s=0}^{2 r_0} (\barq^s 2^j)^{(1 + (2k-1)/M)}\  \ln (\barq^s\, 2^j)  \rkr \\
& = & O \lkr e^{6 (2k-1)\, M \ln M}\   j\  2^{j(1 + (2k-1)/M) }  \rkr, \quad k=1,2,3,4,
\label{final}
\eeqns
proving, thus, \fr{aleph_bounds}. \hfill $\Box$

\bigskip

\noindent{\bf Proof of Lemma \ref{lemm:bound1}.} In what follows, we shall only
construct the proof for the term involving $\bjk$ since the proof for the term involving $\ajk$  is very similar.
Denote 
$$
\Delta_\kappa (j) =
\frac{1}{|C_j|} \sum_{m \in C_j} \lkv \frac{1}{M} \sum_{l=1}^M |g_m(u_l)|^{2} \rkv^{-2 \kappa}\  \lkv \frac{1}{M} \sum_{l=1}^M |g_m(u_l)|^{2 \kappa} \rkv,\ \ \kappa=1,2,
$$
where $\tau_1 (m)$ is given by \fr{eq:box-car-taulower} and \fr{eq:box-car-diophant}.
Note that, by  \fr{alkandblk} and  \fr{coefest}, we have
$$
\hfm - \fm = N^{-1/2} M^{-1}\ \lkv \sumM |\gm(u_l)|^2 \rkv^{-1}\  \lkr \sumM \overline{\gm (u_l)}\, z_{ml} \rkr, 
$$
where $z_{ml}$ are standard Gaussian random variables, independent for different $m$ and $l$. Therefore,
since in the case of Meyer wavelets, $|\psimjk| \leq 2^{-j/2}$ and $|C_j| \asymp 2^j$ 
(see, e.g., Johnstone {\it et al} (2004), p. 565), we derive 
that $\EE |\hbjk - \bjk|^2$ is given by expression \fr{power2}. 
If $\kappa =2$, then 
\begin{align}  
& \EE |\hbjk -
\bjk|^4 =  \, O \lkr \sum_{m \in C_j} \EE |\hfm - \fm|^4 \rkr
+
O \lkr \lkv \sum_{m \in C_j} \EE |\hfm - \fm|^2 \rkv^2 \rkr \nonumber \\
= & O \lkr N^{-2} 2^{-2j} M^{-4} [\tau_1 (m)]^{-4}\  \sum_{m \in C_j} \sum_{l=1}^M |\gm (u_l)|^4 \rkr
+ O \lkr N^{-2} M^{-2} 2^{-2j} [\tau_1 (m)]^{-2} \rkr \nonumber \\
= &  O \lkr 2^{-j} N^{-2} M^{-3} \Delta_2 (j) +
N^{-2} M^{-2} \Delta_1^2 (j) \rkr = O \lkr n^{-2} [M^{-1}  2^{-j} \Delta_2 (j) + \Delta_1^2 (j)] \rkr.
 \label{power4}
\end{align}

Now, recall that $|g_m (u_l)| \asymp |m|^{-1}\, \| m \beta_l\|$  by \fr{sin_dioph}.
Note that, by formula \fr{aleph}, $\aleph_k (j,M)$ is increasing in $r_0$ and recall that, 
by the definition of the Meyer wavelet basis, one has
$|m| \in [(2\pi/3) 2^j,  (8\pi/3) 2^j] \subset  \Omega_j$  with $r_0 = 3 + \log_2 (\pi/3)$ (see  \fr{range_of_m} and \fr{Omj}).
Then, direct calculations yield
\be \label{Deltas}
\Delta_1 (j) = O(2^{j} M \aleph_1(j, M))\ \ \ \mbox{and}\ \ \ \Delta_2 (j) = O(2^{3j} M^4 \aleph_4(j, M)).
\ee
To complete the proof, combine \fr{aleph_bounds}, \fr{power2}, \fr{power4}  and \fr{Deltas} and note that 
$M j^{-1} 2^{-j(1 - 5/M)} = o(1)$ as $n \rightarrow \infty$, since $2^j \geq 2^{j_0}=\ln n$
and $M=M_n \rightarrow \infty$ as $n \rightarrow \infty$.
\hfill $\Box$

\bigskip

\noindent{\bf Proof of Lemma \ref{lemm:bound2}.} 
It is easy to see that $\hbjk -\bjk$ follows a Gaussian distribution with mean zero and variance bounded by 
$\frac{C_{\psi}}{2} (n_M)^{-1}j 2^{j(2+1/M)}$.  Hence,
$$
\PP \lkr |\hbjk - \bjk|^2 \geq \eta^2 (n_M)^{-1} j 2^{j(2+1/M)} \ln n \rkr  \leq  2 \Phi \lkr \frac{\eta}{\sqrt{C_{\psi}}} \sqrt{\ln n} \rkr  
= O \lkr \frac{n^{-\eta^{2}/(2C_{\psi})}}{\sqrt{\ln n}} \rkr,
$$
where $\Phi(\cdot)$ is the cumulative distribution function of a Gaussian random variable with mean zero and variance one. 
\hfill $\Box$

\bigskip

\noindent {\bf Proof of Theorem \ref{th:box-car-slow}.} Due to the orthogonality of the Meyer wavelet basis,
we obtain
$$
\EE \|\hat{f}_n -f\|^2_2= R_0 +R_1+R_2+R_3+R_4,
$$
where
\begin{eqnarray*}
R_0 & =& \sum_{k=0}^{2^{j_0}-1} \EE(\hajk-\ajk)^2, \quad R_1 = \sum_{j=J}^{\infty} \sum_{k=0}^{2^j-1} b_{jk}^2, \\
R_2 &=& \sum_{j=j_0}^{J-1} \sum_{k=0}^{2^j-1} \EE [(\hbjk -\bjk)^2 \II(|\hbjk|\geq \lambda_j)] \,  \II(|\bjk| < \lambda_{j}/2),\\
R_3 &=& \sum_{j=j_0}^{J-1} \sum_{k=0}^{2^j-1} \bjk^2 \PP(|\hbjk| < \lambda_j) \,  \II(|\bjk| \geq 2 \lambda_{j}), \\
R_4 &=& \sum_{j=j_0}^{J-1} \sum_{k=0}^{2^j-1} \EE [(\hbjk -\bjk)^2 \II(|\hbjk| \geq \lambda_j)] \,  \II(|\bjk| \geq \lambda_{j}/2),\\
R_5 &=& \sum_{j=j_0}^{J-1} \sum_{k=0}^{2^j-1} \bjk^2 \PP(|\hbjk| < \lambda_j) \,  \II(|\bjk| < 2 \lambda_{j}).
\end{eqnarray*}
Denote 
$$
\zeta(s,M) = \frac{2(3+1/M)}{2s+3+1/M} 
$$
and observe that $\zeta(s,M) < 2$  for $s >1/\min(p,2)$. 
First, consider the terms $R_0$ and $R_1$. Using Lemma \ref{lemm:bound1}, it is easily seen that
$$
R_0= O \lkr n^{-1} 2^{j_0} \aleph_1 (j_0, M) \rkr = o \lkr n^{-1}\, j_0\, 2^{j_0(2 +1/M)}\ e^{6M\ln M} \rkr =
o((M\, n_M)^{-1}  \ln^3 n) = o \lkr (n_M)^{-\frac{2s}{2s+3+1/M}} \rkr.
$$
Furthermore, it is well-known (see, e.g., Johnstone (2002), Lemma
19.1) that if $f \in \Bpqsa$, then for some positive constant
$c^{*}$, dependent on $p$, $q$, $s$ and $A$ only, we have $
\sumk \bjk^2 \leq c^{*} 2^{-2 j \sstar}$ and,
thus, 
$$
R_1 = O \lkr 2^{-2 J \sstar} \rkr = O \lkr (n_M)^{-2\sstar/(3+1/M)} \rkr.
$$ 
By direct calculations, one can check that if $2 \leq p \leq \infty$, then $\sstar=s$ and hence
$$
R_1 = o \lkr (n_M)^{-2s/(2s+3+1/M)} \rkr.
$$
On the other hand, if $1 \leq p <2$ then $\sstar=s+1/2-1/p$. If 
$\zeta(s,M) \leq p <2$ then $ 2\sstar/(3+1/M) \geq  2s/(2s+3+1/M)$ and, hence,
$$
R_1 = O \lkr (n_M)^{-2s/(2s+3+1/M)} \rkr.
$$ 
Similarly, if $1 \leq p < \zeta(s,M)$, then $ 2\sstar/(3+1/M) \geq  2\sstar/(2\sstar+2+1/M)$ and therefore 
$$
R_1 = O \lkr (n_M)^{-2\sstar/(2\sstar+2+1/M)} \rkr.
$$

Now, consider the term $R_2$. Using Lemma  \ref{lemm:bound1} and Lemma \ref{lemm:bound2}  with $\theta \geq 2$,
formula \fr{ff:newparam}, and the fact that $e^{M \ln M} = o(n^a)$ for any $a>0$ as $n \rightarrow  \infty$,
 after some simple algebra, one derives
\begin{eqnarray*}
R_2 & \leq & \sum_{j=j_0}^{J-1} \sum_{k=0}^{2^j-1} \EE  \lkv (\hbjk -\bjk)^2 \II(|\hbjk - \bjk| \geq  \lambda_j/2) \rkv  
\leq \sum_{j=j_0}^{J-1} \sum_{k=0}^{2^j-1} \sqrt{ \EE [(\hbjk -\bjk)^4} \sqrt{ \PP(|\hbjk - \bjk|^2  \geq \lambda_j^2/4)}  \\
        & =  & O \lkr  \sum_{j=j_0}^{J-1} \sum_{k=0}^{2^j-1} \frac{M e^{21\, M \ln M}\ j 2^{j(2+1/M)}}{n^{1+\theta}} \rkr
        = O \lkr \frac{   2^{J(3+1/M)} \ln n\  e^{15\, M \ln M}}{n_M\,  n^{\theta}} \rkr   = O \lkr (n_M)^{-1} \rkr.
\end{eqnarray*}

For the term $R_3$, again applying  Lemma \ref{lemm:bound2}  with $\theta \geq 2$, obtain
\begin{eqnarray*}
R_3 & \leq & \sum_{j=j_0}^{J-1} \sum_{k=0}^{2^j-1} \bjk^2 \PP(|\hbjk - \bjk| \geq \lambda_j/2)   
          =    o \lkr  \sum_{j=j_0}^{J-1} 2^{-2j\sstar} n^{-\theta}\rkr = o \lkr n^{-1} \rkr.
\end{eqnarray*}

Now, consider the term $R_4$.  Let $j_1$ be such that
$$
2^{j_1}= O \lkr (n_M)^{1/(2s+3+1/M)} (\ln n)^{\xi_0} \rkr,
$$ 
for some real number $\xi_0$.  First, consider the case when $p > \zeta(s,M)$.  Then,
\begin{eqnarray*}
R_4 &\leq & \sum_{j=j_0}^{J-1} \sum_{k=0}^{2^j-1} \EE [(\hbjk -\bjk)^2 \,  \II(|\bjk| \geq \lambda_{j}/2)  
        = R_{41}+R_{42},
\end{eqnarray*}
 where
$$
R_{41} = \sum_{j=j_0}^{j_1} \sum_{k=0}^{2^{j}-1} \EE [(\hbjk -\bjk)^2 \,  \II(|\bjk| \geq \lambda_{j}/2), \quad
 R_{42} =  \sum_{j=j_1+1}^{J-1} \sum_{k=0}^{2^{j}-1} \EE [(\hbjk -\bjk)^2 \,  \II(|\bjk| \geq \lambda_{j}/2).
$$
Then, Lemma \ref{lemm:bound1}  yields
\begin{eqnarray*}
R_{41} &= & O \lkr  \sum_{j=j_0}^{j_1} \sum_{k=0}^{2^j-1} (n_M)^{-1}\ j\ 2^{j(2+1/M)}\rkr
= O \lkr (n_M)^{-2s/(2s+3+1/M)} (\ln n)^{1+\xi_0(3+1/M)}\rkr.
\end{eqnarray*}
For term $R_{42}$,  one derives
\begin{eqnarray*}
R_{42} &= & O \lkr  \sum_{j=j_1+1}^{J-1} \sum_{k=0}^{2^j-1} (n_M)^{-1}\ j\ 2^{j(2+1/M)}\ \frac{|\bjk|^p}{|\lambda_j|^p}\rkr\\    
             &=&  O \lkr  (n_M)^{-(1-p/2)} (\ln n)^{1-p} \sum_{j=j_1+1}^{J-1} 2^{j[(2+1/M)(1-p/2)-\sstar p]} \rkr    
               =  O \lkr (n_M)^{-\rho_1} (\ln n)^{\rho_2}\rkr,    
\end{eqnarray*}
where
$\rho_1=- 2s/(2s+3+1/M)$ and $\rho_2=1-p-\xi_0\, [p \sstar -(2+1/M)(1-p/2)]$. Now, choosing w $\xi_0= -2/(2s+3+1/M)$,
and combining the above terms, one easily arrives at
$$
R_{4}= o \lkr (n_M)^{-\frac{2s}{2s+3+1/M}} (\ln n)^{\frac{2s}{2s+3+1/M}}\rkr.
$$
Now, consider the case when $1 \leq p < \zeta(s,M)$. Note that, same as above,
$$
R_{41}=O \lkr (n_M)^{-2s/(2s+3+1/M)}\ (\ln n)^{1+\xi_0(3+1/M)}\rkr;
$$
but $\xi_0$ does not need to be the   value chosen above. Observe that, since for $R_4$ one has 
$|\bjk| \leq c^{*} 2^{-j \sstar}$ and $|\bjk| > \lambda_j/2$, then, combination of these inequalities requires 
$j \leq j_2$ where  $j_2$ satisfies
$ 
j_2\ 2^{j_2}=O \lkr   n_M/\ln n \rkr^{\frac{1}{2 \sstar+2+1/M}}.
$  
Then,  $|\bjk| \leq \lambda_j/2$ if $j \geq j_2+1$ and 
\begin{eqnarray*}
R_{42}  &=&  O \lkr  (n_M)^{-(1-p/2)}\ (\ln n)^{1-p}\ 2^{j_2\, [(2+1/M)(1-p/2)-\sstar p]} \rkr    
             = O \lkr (n_M)^{\rho_3} (\ln n)^{\rho_4}\rkr,    
\end{eqnarray*}
where  $\rho_3=- 2\sstar/(2\sstar+2+1/M)$ and 
$\rho_4= 2\sstar/(2\sstar+2+1/M) - p/2-[(2+ 1/M)(1- p/2) -p \sstar]$. Noting that, in this case,
$s/(2s+3+1/M) - \sstar/(2\sstar+2+1/M) >0$  and one arrives at $R_{41} = o(R_{42})$ as $n \rightarrow \infty$. Therefore,
$$
R_{4}= O \lkr \lkr \frac{\ln n}{n_M}\rkr^{\frac{2\sstar}{2\sstar+2+1/M}} (\ln n)^{-\frac{p}{2}-[(2+\frac{1}{M})(1-\frac{p}{2}) -p \sstar]}\rkr
= O \lkr \lkr \frac{\ln n}{n_M}\rkr^{\frac{2\sstar}{2\sstar+2+1/M}} \rkr
$$
since the power of $\ln n$ in the expression above is negative.

Finally, consider the term $R_5$.  First, consider the case when $\zeta(s,M) \leq p <2$. Let $j_3$ be such that
$$
2^{j_3}=O \lkr (n_M)^{\frac{s}{\sstar(2s+3+1/M)}} (\ln n)^{\xi_1}\rkr,
$$
for some real number $\xi_1$.  Then,
\begin{eqnarray*}
R_5 &\leq & \sum_{j=j_0}^{J-1} \sum_{k=0}^{2^j-1} \bjk^2 \,  \II(|\bjk| < 2 \lambda_{j}) 
          \leq   R_{51}+R_{52},
\end{eqnarray*}
where
$$
R_{51} = \sum_{j=j_3+1}^{J_1} \sum_{k=0}^{2^{j}-1} \bjk^2 = O \lkr (n_M)^{-\frac{2s}{2s+3+1/M}} (\ln n)^{-2 \sstar \xi_1}\rkr, \quad
R_{52} =  \sum_{j=j_0}^{j_3} \sum_{k=0}^{2^{j}-1} \bjk^2 \,  \II(|\bjk| < 2 \lambda_{j}).
$$
Let 
$$
\Xi(j)=\sum_{k=0}^{2^j-1} \bjk^2\ \II(|\bjk| < 2 \lambda_j).
$$
Note that
$$
\Xi(j)= O \lkr 2^j \lambda_j^2\rkr = O \lkr  j\ 2^{j(3+1/M)} \ln n \ (n_M)^{-1} \rkr
$$
and also
\begin{eqnarray*}
\Xi(j) & = & O \lkr \sum_{k=0}^{2^j-1} |\bjk|^p\ |\bjk|^{2-p}\ \II(|\bjk| < 2 \lambda_j) \rkr  
                   =   O \lkr \lambda_j^{2-p}\ 2^{-jp \sstar} \rkr  \\
                   & = &   O \lkr (n_M)^{p/2-1}\ (\ln n)^{1-p/2}\ j^{1-p/2}\ 2^{j[(2+\frac{1}{M})(1-\frac{p}{2}) -p \sstar]}\rkr.
\end{eqnarray*}
Let $j_4$ be such that
$$
2^{j_4}=O \lkr (n_M)^{\frac{1}{2s+3+1/M}} (\ln n)^{\xi_2}\rkr,
$$
for some real number $\xi_2$. Then
\begin{eqnarray*}
R_{52} &=&\sum_{j=j_0}^{j_4} \Xi(j) +  \sum_{j=j_4+1}^{j_3} \Xi(j)  =  O \lkr \sum_{j=j_0}^{j_4}  j\ 2^{j(3+1/M)}\ \ln n\ (n_M)^{-1} \rkr \\
            & + & O \lkr \sum_{j=j_4+1}^{j_3}  (n_M)^{p/2-1} (\ln n)^{1-p/2}\ j^{1-p/2}\ 2^{j[(2+\frac{1}{M})(1-\frac{p}{2}) -p \sstar]}   \rkr \\
          & =&  O \lkr (n_M)^{-\frac{2s}{2s+3+1/M}} (\ln n)^{2+\xi_2(3+1/M)} \rkr
    +  O \lkr (n_M)^{-\frac{2s}{2s+3+1/M}}\ (\ln n)^{2-p +\xi_2 [(2+\frac{1}{M})(1-\frac{p}{2}) -p \sstar]} \rkr.
\end{eqnarray*}
Since the bound for $R_{52}$ is valid for any value of $\xi_2$, we choose $\xi_2$ which minimizes
$$
\max (2+\xi_2(3+1/M), \,\, 2-p +\xi_2 [(2+1/M)(1-p/2) -p \sstar] ),
$$
i.e., $\xi_2 =  [\sstar +1+1/p +1/(2M)]^{-1}\ (2\sstar+2/p-1).$
Hence
$$
R_{52}= O \lkr (n_M)^{-\frac{2s}{2s+3+1/M}} (\ln n)^{\frac{2\sstar+2/p-1}{\sstar +1+1/p +1/(2M)}} \rkr.
$$
Choose now $\xi_1 =- 2/(2s+3+1/M)$. Then, combining the $R_{51}$ and $R_{52}$ terms, obtain
$$
R_5= O \lkr \lkr \frac{\ln n}{n_M} \rkr^{\frac{2s}{2s+3+1/M}}\ (\ln n)^{\frac{2s}{2s+3+1/M}} \rkr.
$$

Now, consider the case when $1 \leq p < \zeta(s,M)$. Let $j_5$ be such that
$$
2^{j_5} = O \lkr \lkr  \ln n/n_M  \rkr^{\frac{1}{2\sstar+2+1/M}}  (\ln n)^{\xi_3} \rkr,
$$
for some real number $\xi_3$. Then
$$
R_5 \leq R_{51}+ R_{52}+R_{53},
$$
where
$$
R_{51} = \sum_{j=j_5 + 1}^{J-1} \sum_{k=0}^{2^j-1} \bjk^2, \quad
R_{52} =  \sum_{j=j_0}^{j_4} \Xi(j), \quad
R_{53} = \sum_{j=j_4 +1}^{j_5} \Xi(j).
$$
It is immediate that
$$
R_{51}= O \lkr \sum_{j=j_5 + 1}^{J-1} 2^{-2j\sstar} \rkr = 
O \lkr \lkr \frac{\ln n}{n_M} \rkr^{\frac{2\sstar}{2\sstar+2+1/M}}  (\ln n)^{-2 \sstar \xi_3} \rkr.
$$
and that
$$
R_{52} = O \lkr \sum_{j=j_0}^{j_4} \frac{j\ 2^{j(3+1/M)} \ln n}{n_M} \rkr = o \lkr (n_M)^{-\frac{2 \sstar}{2 \sstar+2+1/M}} \rkr. 
$$
After some simple algebra, one obtains
\begin{eqnarray*}
R_{53} & = & O \lkr \sum_{j=j_4 +1}^{j_5}\ (n_M)^{p/2-1}\ (\ln n)^{1-p/2}\ j^{1-p/2}\ 2^{j[(2+\frac{1}{M})(1-\frac{p}{2}) -p \sstar]} \rkr  \\
            &=& O \lkr \lkr \frac{\ln n}{n_M} \rkr^{\frac{2\sstar}{2\sstar+2+1/M}}\  (\ln n)^{1-p/2+\xi_3[(2+\frac{1}{M})(1-\frac{p}{2}) -p \sstar]} \rkr.
\end{eqnarray*}
Choosing
$
\xi_3= - 1/(2 \sstar+2+1/M),
$
and combining the above terms, we arrive at
$$
R_5 = O \lkr \lkr \frac{\ln n}{n_M} \rkr^{\frac{2\sstar}{2\sstar+2+1/M}}\ (\ln n)^{\frac{2\sstar}{2\sstar+2+1/M}} \rkr.
$$
Finally, consider the case, $2 \leq p \leq \infty$. In this case, $j_3=j_4$, and we easily see that
$$
R_5 = O \lkr \lkr \frac{\ln n}{n_M} \rkr^{\frac{2s}{2s+3+1/M}}  \rkr.
$$
Combining all the above expressions, we obtain that, as $n \rightarrow \infty$, 
\be 
R_n (\hat{f}_n) =  \lfi
\begin{array}{ll}
O \lkr \lkr \frac{\ln n}{n_M} \rkr^{\frac{2s}{2s+3+1/M}}\rkr, & {\rm
if}\;\;\; 2 \leq p \leq \infty,
\\ 
O \lkr \lkr \frac{\ln n}{n_M} \rkr^{\frac{2s}{2s+3+1/M}}\rkr (\ln n)^{\frac{2s}{2s+3+1/M}}, & {\rm if}\;\;\; \zeta(s,M) \leq p < 2,
\label{up New} \\
O \lkr \lkr \frac{\ln n}{n_M} \rkr^{\frac{2\sstar}{2\sstar+2+1/M}}\rkr (\ln n)^{\frac{2\sstar}{2\sstar+2+1/M}}, & 
{\rm if}\;\;\; 1 \leq p< \zeta(s,M).
\end{array} \right.
\ee
Now, note that $6/(2s+3) < \zeta(s,M)$ for any $M>0$. Hence, if $p \leq 6/(2s+3)$, then $p < \zeta(s,M)$.
On the other hand, if $p > 6/(2s+3)$, then it is easy to show that for $M$ large enough one has $p > \zeta(s,M)$.
Observe also that $p > 6/(2s+3)$ if and only if $s>3(1/p-1/2)$.

The upper bound in \fr{up New} depends on the choice of $M=M_n$.  Choose $M_n$ of the form \fr{Mnval}.
Then, from the definition of $n_M$ and formulae \fr{A123} and \fr{up New}, it follows  that 
\be  \label{opt_expr}
R_n (\hat{f}_n) = O\big(\exp \big\{ -(A_2+1/M)^{-1} [ A_1 \ln(n - 6M \ln M - \ln M) - A_3 \ln\ln n ] \big\}\big).
\ee
Using Taylor expansion, we write 
$(A_2 + 1/M)^{-1} = A_2^{-1} - M^{-1} A_2^{-2} + M^{-2} A_2^{-3} + O(M^{-3})$ as $M \rightarrow \infty$.
Recalling that $\ln M = \ln \nu + 0.5\, \ln \ln n - 0.5\, \ln \ln \ln n$ and plugging expressions for 
$M$, $\ln M$ and $(A_2 + 1/M)^{-1}$ into the argument of exponent in \fr{opt_expr}, by direct calculations, one derives that 
$$
R_n (\hat{f}_n) = O\big(\exp \lkr - (A_2)^{-1}\, A_1  \ \ln n + \Delta_n \rkr\big),
$$
where, as $n \rightarrow \infty$, 
\begin{eqnarray*}
\Delta_n & = & \sqrt{\ln n} \, \sqrt{\ln \ln n} \lkr \frac{A_1}{A_2} \lkv 3\nu + \frac{1}{A_2 \nu} \rkv \rkr 
-  \frac{3 A_1\nu\ \sqrt{\ln n}\ \ln\ln\ln n}{A_2\ \sqrt{\ln \ln n}} \lkr \frac{2\, \ln \nu}{\ln\ln\ln n} -1 \rkr \\
&-&   \ln \ln n \ \lkr \frac{A_3}{A_2} +  \frac{A_1}{2 A_2} - \frac{3 A_1}{A_2^2} - \frac{A_1}{A_2^3\, \nu^2} \rkr
+  \ln \ln \ln n  \lkr  \frac{A_1}{2 A_2} - \frac{3 A_1}{A_2^2} \rkr + O(1).
\end{eqnarray*}
Now, to complete the proof, note that the main term in $\Delta_n$ is minimized by $\nu = \nu_{opt} = (3  A_2)^{-1/2}$,
and that  $A_2 \geq 2$ for any  $s > 1/\min(p,2)$.
\hfill $\Box$

\section*{Acknowledgments}

Marianna Pensky was supported in part by National Science Foundation
(NSF), grant  DMS-0652524. We would like to thank the Editor and the anonymous referee for useful comments and suggestions on improvements to the presentation of the work.

\end{document}